\documentclass{amsproc}

\usepackage{amscd}
\usepackage{amssymb}
\usepackage{amsthm}
\usepackage{graphicx}
\usepackage{color}
\usepackage[all]{xy}
\usepackage{mathrsfs}
\usepackage{marvosym}
\usepackage{scalerel,amssymb}
\usepackage{hyperref}

\usepackage{tcolorbox}

\usepackage{xcolor}

\usepackage{tikz}
\usetikzlibrary{shapes,arrows}
\usetikzlibrary{positioning}
\usepackage{fixltx2e}
\usetikzlibrary{shapes}

\definecolor{magnolia}{rgb}{0.97, 0.96, 1.0}

\newtheorem{theorem}{Theorem}[section]

\newtheorem{definition}[theorem]{Definition}
\newtheorem{example}[theorem]{Example}

\theoremstyle{remark}
\newtheorem{remark}[theorem]{Remark}
\numberwithin{equation}{section}
\numberwithin{figure}{section}

\newcommand{\NN} {\mathbb{N}}
\newcommand{\ZZ} {\mathbb{Z}}

\newcommand{\RR} {\mathbb{R}}
\newcommand{\CC} {\mathbb{C}}
\newcommand{\PP} {\mathbb{P}}

\newcommand{\cO} {\mathcal{O}}

\newcommand {\Hom}  {\operatorname{Hom}}
\newcommand {\Spec} {\operatorname{Spec}}

\newcommand{\shM} {\mathcal{M}}
\newcommand{\M} {\mathcal{M}}
\newcommand {\lra} {\longrightarrow}
\newcommand\restr[2]{{
  \left.\kern-\nulldelimiterspace 
  #1 
  \vphantom{\big|} 
  \right|_{#2} 
  }}
\newcommand\A{\mathbb A}
\newcommand{\shP} {\mathcal{P}}
\newcommand{\bigslant}[2]{{\raisebox{.2em}{$#1$}\left/\raisebox{-.2em}{$#2$}\right.}}
\newcommand\shO{\mathcal{O}}

\definecolor{White}{rgb}{0.95, 0.95, 0.96}

\numberwithin{equation}{section}

\makeatletter
  \def\title@font{\Large\bfseries}
  \let\ltx@maketitle\@maketitle
  \def\@maketitle{\bgroup%
    \let\ltx@title\@title%
    \def\@title{\resizebox{\textwidth}{!}{%
      \mbox{\title@font\ltx@title}%
    }}%
    \ltx@maketitle%
  \egroup}
\makeatother

\begin{document}

\title[Quiver DT/log Gromov--Witten theory]{Quiver DT Invariants and Log Gromov--Witten Theory of Toric Varieties}



\author[H.\,Arg\"uz]{H\"ulya Arg\"uz}
\address{University of Georgia, Department of Mathematics, Athens, GA 30605}
\email{Hulya.Arguz@uga.edu}

\subjclass[2000]{Primary }

\date{}

\begin{abstract}
We review how log Gromov--Witten invariants of toric varieties can be used to express quiver Donaldson--Thomas invariants in terms of the simpler attractor Donaldson--Thomas invariants.
This is an exposition of joint work with Pierrick Bousseau \cite{AB,AB2}. 
\end{abstract}

\maketitle

	\setcounter{tocdepth}{1}
	\tableofcontents
	\setcounter{section}{-1}
	\section{Introduction}
Donaldson-Thomas (DT) invariants of a quiver with potential can be expressed in terms of simpler attractor DT invariants by a universal formula. The coefficients in this formula are calculated combinatorially using attractor flow trees. 
In this paper, we give an introduction to the main result of \cite{AB2} showing that these coefficients are genus 0 log Gromov--Witten invariants of $d$-dimensional toric varieties, where $d$ is the number of vertices of the quiver.



\section{Quiver DT invariants}

\subsection{Quiver representations and moduli spaces}

A \textbf{quiver} is a finite oriented graph $Q$. We denote by
\begin{itemize}
\item[] $Q_0$: the finite set of vertices of $Q$,
\item[]  $Q_1$: the finite set of arrows (oriented edges) of $Q$,
\item[] $s: Q_1 \to Q_0$ the source map which maps an arrow to its source, and;
\item[] $t: Q_1 \to Q_0$ the target map, which maps an arrow to its target.
\end{itemize}

We illustrate a quiver in Figure \ref{Fig:1}. 

\begin{figure}[ht]
\center{\scalebox{.45}{\input{Quiver.pspdftex}}}
\caption{A quiver}
\label{Fig:1}
\end{figure}

\begin{definition}
A \textbf{representation} of a quiver $Q$ is an assignment of
\begin{itemize}
\item[(i)] a finite-dimensional vector space $V_i$ over $\CC$ for each vertex $i \in Q_0$, and
\item[(ii)] a $\CC$-linear transformation $f_{\alpha} \in \mathrm{Hom}(V_{s(\alpha)},V_{t(\alpha)})$ for each arrow $\alpha\in Q_1$.
\end{itemize}
\end{definition}

\begin{figure}[ht!]
\center{\scalebox{.3}{\input{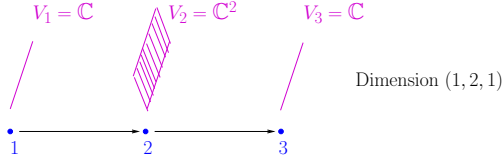}}}
\caption{A quiver with dimension vector $\gamma=(1,2,1)$}
\label{Fig:2}
\end{figure}

Representations of $Q$ naturally form an abelian category: a morphism between $(V_i,f_{\alpha})$ and $(V_i', f_{\alpha}')$ is the data of linear maps $V_i \rightarrow V_i'$ intertwining the maps $f_{\alpha}$ and $f_{\alpha}'$. In what follows, we fix
\[N=\ZZ^{Q_0}=\bigoplus_{i \in Q_0} \ZZ e_i \,.\] 
The \textbf{dimension} of a quiver representation $(V_i, f_{\alpha})$ is defined as the vector 
\[\gamma = (\gamma_i)_{i \in Q_0} \in N\,,\] 
where $\gamma_i :=\dim V_i$ -- see Figure \ref{Fig:2} for an illustration.

There is a notion of stability due to King in the context of quiver representations, defined as follows \cite{MR1315461}. Denote $M: = \mathrm{Hom}(N,\ZZ)$ the dual lattice to $N$ and set 
\[M_{\RR} =  \mathrm{Hom}(N,\RR) = M \otimes \RR.\]

\begin{definition}[King's stability]
A \textbf{stability parameter} for $\gamma \in N$ is a point $\theta \in \gamma^{\perp}:=\{ \theta \in M_\RR\,, \theta(\gamma)=0\} \subset M_{\RR}$. A quiver representation $V$ of
dimension $\gamma$ is $\theta$-\textbf{stable} (resp.\ $\theta$-\textbf{semistable}) if for all 
non-zero strict subrepresentation $V'$ of $V$ we have $\theta(\mathrm{dim}(V')) < 0 $
(resp.\ $\theta(\mathrm{dim}(V')) \leq 0 $).
\end{definition}

Given a dimension vector $\gamma \in N$ and a stability parameter $\theta \in M_\RR$, one obtains by geometric representation theory \cite{MR1315461} a moduli space
$\mathcal{M}_{\gamma}^{\theta}$ of S-equivalence classes of $\theta$-semistable quiver representations of $Q$ dimension $\gamma$. It is a quasi-projective variety over $\CC$, which is actually projective if $Q$ is acyclic.

When working with a non-acyclic quiver $Q$, we consider the additional data of a \textbf{potential} $W$, that is, a formal linear combination of oriented cycles of $Q$. A choice of potential $W$ defines a \textbf{trace function} $ \mathrm{Tr}(W)_{\gamma}^{\theta}:  \mathcal{M}_{\gamma}^{\theta} \to \CC $ on each of the moduli spaces 
$\mathcal{M}_{\gamma}^{\theta}$: if $c$ is an oriented cycle of arrows $\alpha_1,\dots,\alpha_n$ of $Q$, 
$ \mathrm{Tr}(c)_{\gamma}^{\theta}:  \mathcal{M}_{\gamma}^{\theta} \to \CC $ is defined by
\[  V = (V_i, f_{\alpha}) \longmapsto \mathrm{Tr}(f_{\alpha_n} \circ  \ldots \circ f_{\alpha_1} ) \]
and we set
$$\mathrm{Tr}(W)_{\gamma}^{\theta}  = \sum_{c} \lambda_c \mathrm{Tr}(c)_{\gamma}^{\theta}$$
if $W$ is the linear combination $\sum_c \lambda_c c$ of oriented cycles.

\subsection{DT invariants}
Let $(Q,W)$ be a quiver with potential, $\gamma \in N$
a dimension vector, and $\theta$ a generic stability parameter for $\gamma$. The DT (Donaldson-Thomas) invariant $\Omega_\gamma^{+,\theta}$ is an integer which is a virtual count of the critical points of the trace function $\mathrm{Tr} (W)_\gamma^\theta$ on the moduli space $\mathcal{M}_{\gamma}^{\theta}$ of $\theta$-semistable representations of dimension $\gamma$.

Let $(Q,W)$ be a quiver with potential, $\gamma \in N$
a dimension vector, and $\theta$ a generic stability parameter for $\gamma$. The \textbf{DT invariant} $\Omega_\gamma^{+,\theta}$ is defined as follows.
If the $\theta$-stable locus in $\mathcal{M}_{\gamma}^{\theta}$ is empty, we have $\Omega_\gamma^{+,\theta}=0$ , and else,  $\Omega_\gamma^{+,\theta}$ is the Euler characteristic of $\mathcal{M}_{\gamma}^{\theta}$ valued in the perverse sheaf obtained by applying the vanishing-cycle functor for $\mathrm{Tr} (W)_\gamma^\theta$ to the intersection cohomology sheaf \cite{davison2015donaldson, MR4132957, MR4000572}:
\[ \Omega_\gamma^{+,\theta}
:=e(M_\gamma^\theta, \phi_{\mathrm{Tr}(W)_\gamma^\theta}(IC[-\dim \mathcal{M}_{\gamma}^{\theta}]))\,.\]

For example, if $\gamma$ is primitive and $W=0$, then $\Omega_\gamma^{+,\theta}$ is simply the topological Euler characteristic\footnote{Here, we are considering the ``unsigned" DT invariants $\Omega_\gamma^{+,\theta}$, and not the signed DT invariants $\Omega_\gamma^{\theta}$, which are equal to $(-1)^{\dim \mathcal{M}_{\gamma}^{\theta}}e(\mathcal{M}_{\gamma}^{\theta})$ when $\gamma$ is primitive and $W=0$.}
of $\mathcal{M}_{\gamma}^{\theta}$.
The DT invariants 
$\Omega_\gamma^{+,\theta}$ can be repackaged into the \textbf{rational DT invariants} 
\[
\overline{\Omega}_\gamma^{+,\theta} := \sum_{\substack{\gamma' \in N\\ 
    \gamma=k \gamma',\, k\in \ZZ_{\geq 1}}} \frac{(-1)^{k-1}}{k^2} \Omega_{\gamma'}^{+,\theta} \,,
\]
 which can also be defined using the motivic Hall algebra \cite{JoyceSong, kontsevich2008stability, MR2650811,MR2801406}.



\begin{figure}[ht!]
\center{\scalebox{.50}{\input{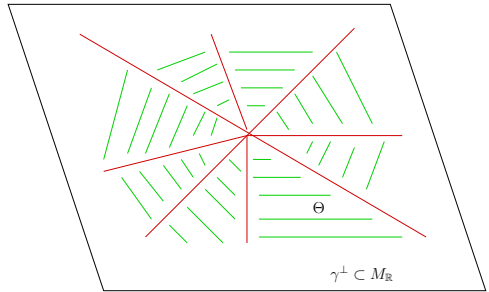}}}
\caption{A wall and chamber structure}
\label{Fig:3}
\end{figure}

The quiver DT invariants $\Omega_\gamma^{\theta,+}$ have a piecewise constant dependence on $\theta \in \gamma^{\perp}$: there is a natural wall and chamber structure on $M_\RR$ as illustrated in Figure \ref{Fig:3}, and while the quiver DT invariants remain constant as long as $\theta$ remains inside a chamber, they change by a universal \textbf{wall-crossing formula} \cite{JoyceSong, kontsevich2008stability} once we cross a wall and pass to another chamber -- see \ref{Fig:3}.


\begin{example} \label{ex_n_K}
Let $Q$ be the $n$-Kronecker quiver, consisting of two vertices connected by $n$ arrows as illustrated in Figure \ref{Fig:4}. Fix the dimension vector $\gamma := \mathrm{dim}(V) = (1,1) \in N$ and a choose a stability parameter $\theta = (\theta_1,-\theta_1) \in \gamma^{\perp} \subset M_{\RR}$.
The data of $n$ linear maps $\CC \to \CC$ is given by a tuple of $n$ complex numbers $(\xi_1,\ldots, \xi_n)$, and a global rescaling of $(\xi_1,\ldots, \xi_n)$ by a non-zero constant produces isomorphic representations. For $\theta_1 > 0$, $\theta$-semistability requires $(\xi_1,\ldots, \xi_n) \neq (0,\ldots,0)$, and in this case, the moduli space of $\theta$-semistable quiver respresentations $\mathcal{M}_{\gamma}^{\theta}$ is $\CC\PP^{n-1}$, and so \[\Omega_\gamma^{\theta,+}=e(\CC\PP^{n-1})=n.\] On the other hand, for  $\theta_1 < 0$ we have $\mathcal{M}_{\gamma}^{\theta} = \emptyset $, and so $\Omega_\gamma^{\theta,+}=0$. 


\begin{figure}[ht!]
\center{\scalebox{.50}{\input{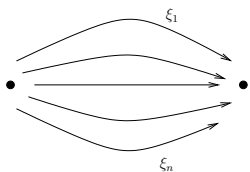}}}
\caption{The $n$-Kronecker quiver}
\label{Fig:4}
\end{figure}

\end{example}

Quiver DT invariants are extensively studied both in the physics as well as in the mathematics literature. In physics they appear in the context of supersymmetric quantum mechanics, and of BPS states in $\mathcal{N}=2$ four-dimensional quantum field theories and string theory compactifications \cite{MR3268234, MR2567952, MR3036499}, whereas in mathematics they are closely related to geometric DT theory, which concerns counts of coherent sheaves or special Lagrangians in Calabi--Yau threefolds \cite{JoyceSong, kontsevich2008stability, MR1818182} -- see Table \ref{Tab:1}.

\begin{figure}[ht!]
\center{\scalebox{.55}{\input{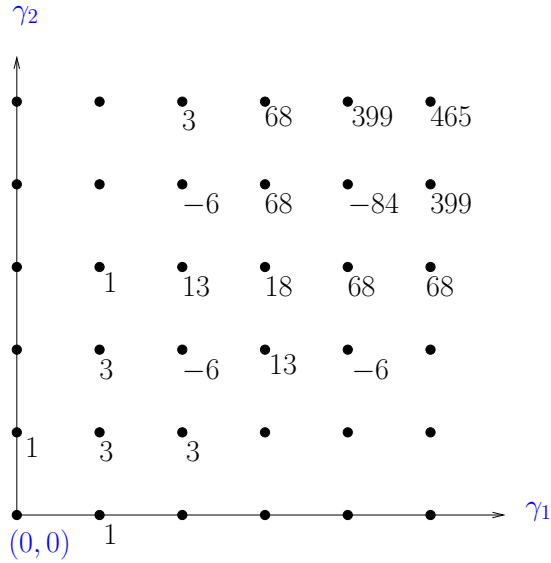}}}
\caption{Values of $\Omega_{\gamma}^{\theta}$ for the $3$-Kronecker quiver}
\label{Fig:BPS}
\end{figure}

\begin{example}
The $3$-Kronecker quiver appears in the description of the BPS spectrum of $\mathcal{N}=2$, four-dimensional $SU(3)$ super Yang-Mills theory \cite{galakhov2013wild} In this context, explicit values of quiver DT invariants are calculated --see \cite[A.2]{galakhov2013wild} and Figure \ref{Fig:BPS}.
\end{example}

\begin{flushleft}
\begin{table}
\begin{tikzpicture}[
roundnode/.style={circle, draw=green!60, fill=green!5, very thick, minimum size=7mm},
squarednode/.style={rectangle, draw=purple!60!black, fill=magnolia!5!White, very thick, minimum size=5mm},
]

\noindent

\node[squarednode,rounded corners=5pt,node distance=0.5cm]   (coh)       {\# Coherent sheaves in CY$3$'s};

\node[squarednode,rounded corners=5pt,node distance=0.5cm]  [ right =of coh]    (Lag)       {\# Special Lagrangians in CY$3$'s};

\node[squarednode,rounded corners=5pt,node distance=0.5cm]  [ below right = 1cm and -1.5cm  of coh]       (geometric)       {Geometric DT theory};

\node[squarednode,rounded corners=5pt,node distance=0.8cm]      (refined)       [ below =of geometric] {DT invariants of quiver representations};

\node[squarednode,rounded corners=5pt,node distance=0.8cm]      (super)       [ below =of refined] {Supersymmetric quantum mechanics};

\node[squarednode,rounded corners=5pt,node distance=0.8cm]      (ground)       [ below left = 1cm and -3cm of super] {Supersymmetric ground states};

\node[squarednode,rounded corners=5pt,node distance=0.8cm]      (BPS)       [ right = of ground] {\# BPS particles/black holes};

Lines
\draw[->] (refined.north) -- (geometric.south);
\draw[->] (refined.south) -- (super.north);
\draw[->] (geometric.north) -- (coh.south);
\draw[->] (geometric.north) -- (Lag.south);
\draw[->] (super.south) -- (ground.north);
\draw[->] (super.south) -- (BPS.north);
\end{tikzpicture}
\caption{DT invariants of quivers in maths and physics}
\label{Tab:1}
\end{table}
\end{flushleft}


\subsection{The flow tree formula}
\label{sec_flow_tree}
Quiver DT invariants can be expressed in terms of a simpler set of quiver DT invariants, known as attractor DT invariants defined as follows.
Let $\langle -,-\rangle$ be the skew-symmetrized Euler form on \[N=\ZZ^{Q_0}=\bigoplus_{i\in Q_0} \ZZ e_i .\]
It is the skew-symmetric form on $N$ defined by $\langle e_i, e_j\rangle=a_{ij}-a_{ji}$, where $a_{ij}$ is the number of arrows from $i$ to $j$ in $Q$.
Given $\gamma \in N$, $\langle \gamma,-\rangle$ is a point in the hyperplane $\gamma^\perp$ of stability parameters for $\gamma$, referred to as the \textbf{attractor point}. 


\begin{definition}
\cite{AlexandrovPioline, MR3330788, mozgovoy2020attractor}  
Fix a dimension vector $\gamma \in N$. The \textbf{attractor DT invariant} for $\gamma$ is defined by $\Omega_\gamma^\star := \Omega_\gamma^\theta$, where $\theta$ is a small generic perturbation in $\gamma^\perp$ of the attractor point $\langle \gamma, -\rangle$ \footnote{One can show that the attractor DT invariant is independent of the choice of perturbation, see \cite[Thm 3.7]{mozgovoy2020attractor}.}.
\end{definition}

The attractor DT invariants turn out to be particularly simple in many situations. 

\begin{example}
It is shown by Bridgeland \cite{bridgeland}  that if $Q$ is acyclic then 

\[ \Omega_\gamma^\star= \begin{cases} 
      1 & \mathrm{if} \,\ \gamma = e_i \,\,\text{for some}\,\, i\in Q_0 \\
      0 & \mathrm{otherwise}
   \end{cases}
\]
This theorem has been generalized by Lang Mou to some situations when $Q$ is not acyclic, but admits a so called red to green sequence \cite{mou2019scattering}.
\end{example}

\begin{example}
Given a toric Calabi-Yau 3-fold, one can construct a quiver with a potential $(Q,W)$ such that
\[D^bRep(Q,W) \cong D^bCoh(X),\]
where $D^bRep(Q,W)$ is the bounded derived category of representations of $(Q,W)$ and $D^bCoh(X)$ is the bounded derived category of coherent sheaves on $X$ \cite{mozgovoy2009crepant}. 

Beaujard--Manschot--Pioline \cite{beaujard2020vafa} and Mozgovoy--Pioline \cite{mozgovoy2020attractor} conjectured in this situation that one gets $\Omega_\gamma^\star= 0$ unless either $\gamma=e_i$ for some $i \in Q_0$ (in which case one has $ \Omega_\gamma^\star = 1$) or $\gamma$ is in the kernel of $\langle -,-\rangle$\footnote{A more precise conjecture including the values of all non-trivial attractor inariants has been formulated by Descombes \cite{MR4524190}.}. This conjecture is proven for local $\PP^2$ by Bousseau--Descombes--Le Floch--Pioline \cite{bousseau2022bps}. 

More generally, attractor DT invariants are expected to take a simple form in physics set-ups without dynamical gravity, such as string compactifications on non-compact Calabi-Yau 3-folds. For theories with dynamical gravity, such as string compactifications on compact Calabi-Yau 3-folds, one typically expects attractor invariants to be sufficiently non-trivial to account for the entropy of black holes \cite{MR3036499}.
\end{example}

\begin{figure}[ht]
\center{\scalebox{.5}{\input{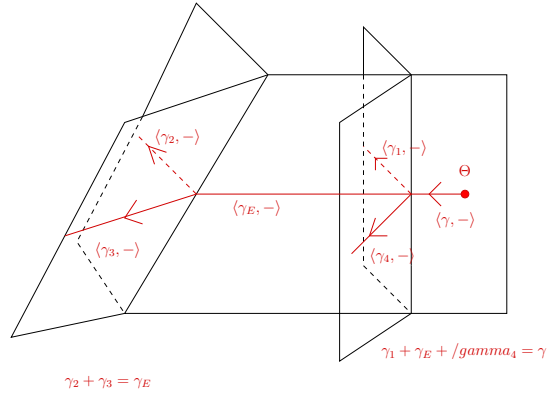}}}
\caption{An attractor tree with higher valent vertices}
\label{Fig:attractor KS}
\end{figure}

Fix a dimension vector $\gamma \in N$, and a generic stability parameter $\theta \in \gamma^{\perp}$. By an iterative application of the wall-crossing formula \cite{MR3330788}, quiver DT invariants can be obtained in terms of attractor DT invariants by a formula of the form
\begin{equation}
\label{Eq: KS attractor}
\overline{\Omega}_\gamma^\theta=\sum_{\gamma=\gamma_1+\dots +\gamma_r}
\frac{1}{|\mathrm{Aut}((\gamma_i)_i)|} 
F_r^\theta(\gamma_1,\dots,\gamma_r) \prod_{i=1}^r
\overline{\Omega}_{\gamma_i}^{\star}\,,    
\end{equation}
where  $|\mathrm{Aut}((\gamma_i)_i)|$ is the order of the group of permutation symmetries of the decomposition 
$\gamma=\gamma_1+\dots+\gamma_r$, and  the coefficients  $F_r^\theta(\gamma_1,\dots,\gamma_r)$ are sums of contributions from attractor trees with leaves decorated by $\gamma_1,\dots,\gamma_r$ and with root at $\theta$:
\[F_r^\theta(\gamma_1,\dots,\gamma_r) =\sum_{T \in \mathcal{T}^\theta_{\gamma_1,\dots,\gamma_r}} F_{r,T}^\theta (\gamma_1,\dots,\gamma_r) .\]
These trees, as illustrated in Figure \ref{Fig:attractor KS}, are oriented trees with one root vertex corresponding to the point $\theta \in M_{\RR}$ and edges have direction vectors determined by the decompositions of the dimension vector $\gamma$.

By construction, attractor trees satisfy the ``tropical balancing condition'' -- that is, weighted directions around edges add up to zero. Hence, they can be viewed as tropical trees. As illustrated in Figure \ref{Fig:attractor KS}, these trees have in general vertices of valency bigger than $3$, 
and calculating their contributions in this case is technically challenging, as it requires the application of the wall-crossing formula.

\begin{figure}[ht]
\center{\scalebox{.5}{\input{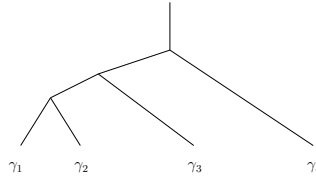}}}
\caption{A binary tree with leaves labelled by $\gamma_1,\ldots,\gamma_n$, for $n=4$}
\label{trees}
\end{figure}

In joint work with Bousseau \cite{AB}, we proved the \textbf{flow tree formula} conjectured by Alexandrov-Pioline \cite{AlexandrovPioline} 
computing the coefficients $F_r^\theta(\gamma_1,\dots,\gamma_t)$
in terms of binary trees as illustrated in Figure \ref{trees}, counted with simple tropical multiplicities attached to the trivalent vertices. These binary trees are obtained as generic trivalent perturbations of the non-trivalent attractor trees. 



In the following section we explain how to interpret the coefficients $F_r^\theta(\gamma_1,\dots,\gamma_r)$ in \eqref{Eq: KS attractor} geometrically in terms of counts of rational curves in toric varieties, by relating contributions of binary trees to such curve counts.

\section{Log Gromov--Witten invariants of toric varieties}
The main result of \cite{AB} relates DT invariants of quivers to log Gromov--Witten invariants given by counts of log curves toric varieties. After briefly explaining the toric varieties we consider and the enumerative problem of curves, we provide further details of log geometry and log Gromov--Witten theory. In the final section we state the relationship between counts of log curves in toric varieties and quiver DT invariants, and provide a brief summary of the proof.

\subsection{From quivers to enumerative geometry in toric varieties}
\label{sec_enum}
Fix a quiver $Q$, a dimension vector $\gamma \in N$ and a decomposition \[\gamma = \gamma_1 + \ldots + \gamma_r.\] 
Recall that we denote by $\langle -,-\rangle$ the skew-symmetrized Euler form on $N$. We assume that $\langle \gamma_i, -\rangle \neq 0$ for all $1\leq i\leq r$, and $\langle \gamma,-\rangle \neq 0$.

Let $\Sigma$ be a fan in $M_\RR$ of a smooth projective toric variety $X_\Sigma$ containing the rays $\RR_{\geq 0}\langle \gamma_i, -\rangle$ for all $1 \leq i\leq r$. We consider hypersurfaces $H_i \subset D_i$ defined by equations 
\[(z^{\frac{\gamma_i}{|\gamma_i|}}+c_i)|_{D_i}=0\, ,\] where $|\gamma_i|$ is the divisibility of $\gamma_i$ in $N$, and $c_i$ is a general constant.

\begin{figure}[ht]
\center{\scalebox{.45}{\input{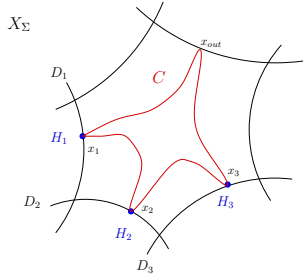}}}
\caption{A log curve in a toric variety with prescribed tangencies along the toric boundary}
\label{curves}
\end{figure}

The count of log curves which we show are related to quiver DT invariants in \cite{AB} are given by appropriately defined counts of marked genus $0$ stable maps 
\[f: (C, \{ p_1, \ldots , p_{r+1} \} ) \to X_{\Sigma}\] satisfying
\begin{itemize}
    \item[(i)] $f(p_i) \in H_i$ for all $1 \leq i \leq r$ 
  \item[(ii)] The contact order of $f$ with $D_i$ at $p_i$ is the divisibility of $\langle \gamma_i , - \rangle $ in $M$.
  \end{itemize}
See Figure \ref{curves} for an illustration.

To count such curves with prescribed tangency conditions, we use log Gromov--Witten theory developed by Abramovich--Chen \cite{AbramovichChen} and Gross--Siebert \cite{GSlogGW}. In the following sections we briefly review log geometry and log Gromov--Witten theory before stating a precise relationship between counts of log curves and quiver DT invariants.

\subsection{Log schemes}
We start reviewing a couple of definitions from log geometry to fix our notation. For familiarity with log geometry we refer to \cite{MR1463703,O}. All the monoids considered below are assumed to be commutative.

\begin{definition}
\label{Def: log structure}
Let $X$ be a scheme. A \textbf{prelog structure} on $X$ is a sheaf of monoids $\shM$ on $X$ together with a homomorphism of monoids 
$\alpha \colon \shM \lra \mathcal{O}_X$ where we consider the structure sheaf $\mathcal{O}_X$ as a monoid with respect to multiplication. 
A pre log structure on $X$ is called a \textbf{log structure} if $\alpha$ induces an isomorphism
\[
\restr{\alpha}{\alpha^{-1}(\mathcal{O}_X^{\times})}:{\alpha}^{-1}(\mathcal{O}_X^{\times})\lra \mathcal{O}_X^{\times}.
\]
We call a scheme $X$ endowed with a log structure a \textbf{log scheme}, and denote the log structure on $X$ by $(\mathcal{M}_X, \alpha_X)$, 
or sometimes by omitting the structure homomorphism from the notation, just by $\mathcal{M}_X$. 
We denote a scheme $X$, endowed with a log structure by $(X, \mathcal{M}_X)$.
\end{definition}

All log schemes below are assumed to be fine and saturated (fs) as in
\cite[I,\S1.3]{O}.
The combinatorial information of the log structure is encoded via tropical geometry using the ghost sheaf as defined below.

\begin{definition}
Let $X$ be a log scheme. The \textbf{ghost sheaf} of $X$ is the sheaf on $X$ defined by
\[
\overline{\M}_{X} \colon = \M_X /{\alpha^{-1}(\mathcal{O}_{X}^{\times})}\,.
\]
The \textbf{tropicalization} $\Sigma(X)$ of $X$ is an abstract polyhedral
cone complex:
$\Sigma(X)$ is the collection of cones  $\sigma_{\bar x}:=
\Hom(\overline{\shM}_{X,\bar x},\RR_{\ge 0})$ for every 
geometric point $\bar x\rightarrow X$,
glued together using the natural generization and specialization maps,
 see \cite[\S2.1]{ACGSI} for details.
\end{definition}

\begin{example} 
\label{The divisorial log structure}
Let $D\subset X$ be a divisor. Let $j \colon X\setminus D\hookrightarrow X$ denote the inclusion map. 
The \emph{divisorial log structure} on $X$ is the pair $(\shM_{(X,D)},\alpha_X)$, 
where $\shM_{(X,D)}$ is  the sheaf of regular functions on $X$, that restrict to units on $X\setminus D$. 
That is,
\[\shM_{(X,D)} = j_*(\mathcal{O}_{X\setminus D}^\times)\cap \mathcal{O}_X\] 
The structure homomorphism $\alpha_X$ is given by the inclusion $\alpha_X:\shM_{(X,D)} \hookrightarrow \mathcal{O}_X$.
\end{example}
\begin{remark}
The log structure we consider on the toric variety $X_{\Sigma}$ constructed from the data of a quiver, as discussed in \ref{sec_enum}, is the divisorial log structure defined by the toric boundary divisor $D_{\Sigma}$ in $X_{\Sigma}$, that is, the natural anti-canonical divisor formed by the union of divisors fixed under the torus action on $X_{\Sigma}$.

Note that, for any toric variety $X_{\Sigma}$ endowed with the divisorial log structure defined by the toric boundary divisor, the tropicalization of the resulting log scheme is naturally identified with the toric fan $\Sigma$ associated to $X$.
\end{remark}

Analogous to assigning a sheaf to a presheaf, we can assign a log structure to a prelog structure, using a fibered coproduct, as follows.
\begin{definition}
\label{Def: Log pre log}
Let $\alpha \colon \shP\lra \mathcal{O}_X$ be a prelog structure on $X$.  
We define the \emph{log structure associated to the prelog structure} $(\shP,\alpha)$ on $X$ as follows. Set 
\[\shP^a = \bigslant{\shP\oplus \mathcal{O}_X^{\times}}
{\{(p,\alpha(p)^{-1})\,\big|\, p\in \alpha^{-1}(\mathcal{O}^\times_X)\big\}}\]
and define the structure homomorphism 
$\alpha^a \colon  \shP^a  \longrightarrow  \mathcal{O}_X$ by
\begin{eqnarray}
\nonumber
\alpha^a(p,h) & = h\cdot \alpha(p)
\end{eqnarray}
One can easily check that $(\shP^a,\alpha^a)$ is a log structure on $X$.
\end{definition}

Let $Y$ be a log scheme, and let $f\colon X\to Y$ be a scheme theoretic morphism. Then, the log structure on $Y$, given by $\alpha_Y\colon \shM_Y\lra \mathcal{O}_Y$,
induces a log structure on $X$ defined as follows. First define a prelog structure,
by considering the composition
\[  f^{-1}\shM_Y\to f^{-1}\mathcal{O}_Y \to \mathcal{O}_X \,.\]
Then, we endow $X$ with the log structure associated to this prelog structure, as in Definition \ref{Def: Log pre log}. We refer to this log structure as the \emph{pull back log structure} or the \emph{induced log structure} on $X$, and denote it by $\shM_X = f^* \shM_Y$.

\begin{example}
\label{Ex: trivial log point}
For every scheme $X$, $\mathcal{M}_X = \cO_X^{\times}$ define a log structure on $X$, called the \emph{trivial log structure}. In particular, taking $X= \Spec \CC$, we call 
$(\Spec \CC,\CC^{\times})$ the \emph{trivial log point}.
\end{example}

\begin{example}
\label{Ex: standard log point}
Let $X = \Spec \CC$. Define $\mathcal{M}_X = \CC^{\times} \oplus \NN$, and $\alpha_X\colon \M_X \to \CC$ as follows.
\begin{center}
$\alpha_X ( x , n ) = \left\{
	\begin{array}{ll}
		x  & \mbox{if } n = 0 \\
		0 & \mbox{if } n \neq 0
	\end{array}
\right.
$
\end{center}
The corresponding log scheme $O_0 = (\Spec \CC, \CC^{\times} \oplus \NN)$ is called the \emph{standard log point}. 
One can check that the log structure on the standard log point is the same as the pull-back of the divisorial log  structure 
on $\A^1$ with the divisor $D=\{0\}\subset \A^1$. One also obtains a more general log point $(\Spec \CC, \CC^{\times} \oplus Q)$ by replacing $\NN$ by a more general monoid $Q$ such that $Q^{\times}=\{0\}$. It is the pullback to the torus fixed point of the divisorial log structure on the affine toric variety $\Spec \CC[Q]$. Its tropicalization is the cone $\Hom(Q,\RR_{\geq 0})$.
\end{example}

\begin{definition}
\label{Def: log morphism}
A \textbf{morphism of log schemes} $f \colon (X,\mathcal{M}_X) \to (Y,\mathcal{M}_Y)$ is a morphism of schemes $f \colon X \to Y$
along with a homomorphism of sheaves of monoids $f^\sharp \colon f^{-1}\mathcal{M}_Y \to \mathcal{M}_X$ such that
the diagram
\[\begin{CD}
f^{-1}\M_Y @>{f^\sharp}>> \M_X\\
@V{\alpha_Y}VV@VV{\alpha_X}V\\
f^{-1}\shO_Y@>{f^*}>>\shO_X.
\end{CD}\]
is commutative. Here, $f^*$
is the usual pull-back of regular functions defined by the
morphism $f$. Given a morphism of log spaces $f: (X,\M_X) \to (Y,\M_Y)$, we denote by $\underline{f}\colon X \to Y$ 
the underlying morphism of schemes. By abuse of notation, the underlying morphism on topological spaces is also denoted by $\underline{f}$. 
\end{definition}


We refer to \cite[Defn. $3.23$]{Gross} for the notion of \emph{log smooth morphism}.
We recall that a toric variety is log smooth over the trivial log point. 
If $\pi \colon \mathcal{X} \rightarrow \A^1$ is a toric degeneration of toric varieties, then 
$\pi$ can be naturally viewed as a log smooth morphism. In restriction to $t=0$, we get a log 
smooth morphism $X_0 \rightarrow O_0$, that is, $X_0$ is log smooth over the standard log 
point $O_0$.
We also refer to \cite[Remark 3.25]{Gross} for the notion of 
\textbf{integral} morphism of log schemes.

\subsection{Stable log maps}
\label{Sec: Log curves}

We recall that a $\ell$-marked stable map with target a scheme $X$ is 
a map $f \colon C \rightarrow X$, where $C$ is a proper nodal curve
with $\ell$ marked smooth points $\mathbf{x}=(x_1,\ldots,x_\ell)$,
such that the group of automorphisms of $C$ fixing $\mathbf{x}$
and commuting with $f$ is finite.

The notion of stable map has a natural generalization to the setting of logarithmic 
geometry \cite{AbramovichChen, GSlogGW}. 

\begin{definition}
\label{logCurve}
Let $X \rightarrow B$ be a log morphism between log schemes.
A \textbf{$\ell$-marked stable log map with target $X \rightarrow B$}
is a commutative diagram of log schemes
\[\begin{CD}
C @>{f}>> X\\
@V{\pi}VV@VV{}V\\
W@>{}>>B\,,
\end{CD}\]
together with a tuple of sections $\mathbf{x} = (x_1,\dots,x_{\ell})$
of $\pi$,
where
$\pi$ is a proper log smooth and integral morphism of
log schemes, such
that, for every geometric point $w$ of $W$, the restriction of $\underline{f}$ to $w$
with the marked points $\underline{\mathbf{x}}(w)$ is an ordinary stable map, and such that, if $U \subset C$ is the non-critical locus of $\underline{\pi}$, we have 
$\overline{\shM}_C|_U \simeq \underline{\pi}^{*} \overline{\shM}_W
\oplus \bigoplus_{j=1}^{\ell} (x_j)_{*}\NN_W$.
\end{definition}

When the scheme underlying $W$ in Definition \ref{logCurve}
is a point, then $W$ is a log point $(\Spec \CC, \CC^{*} \oplus Q)$ as in Example \ref{Ex: standard log point}.
The general enumerative theory of stable log maps, called log Gromov--Witten theory and developed in \cite{GSlogGW, AbramovichChen}, is based on the notion of \textbf{basic stable log map}.
For a basic stable log map over a log point 
$(\Spec \CC, \CC^{*} \oplus Q)$, the monoid $Q$ is the so-called 
\textbf{basic monoid}. The basic monoid is uniquely determined by the combinatorial type of the stable log map and has a natural tropical geometric interpretation: the dual cone $\Hom(Q,\RR_{\geq 0})$
is the base of a universal family of tropical curves with given combinatorial type. 
We provide a more precise description of the combinatprial type of a stable log map, and of basic monoids below.

The \textbf{combinatorial type} $\tau$ of a stable log map 
$f \colon C/W \rightarrow X/S$
consists of:
\begin{itemize}
\item[(i)]
The dual intersection graph $G=G_C$ of $C$, with set of vertices $V(G)$, set of edges $E(G)$, and set of legs $L(G)$.
\item[(ii)]
The map $\sigma: V(G)\cup E(G)\cup L(G)\rightarrow \Sigma(X)$ mapping $x\in C$ to
$\big(\overline{\M}_{X,f(x)}\big)_\RR^\vee$.
\item[(iii)]
The contact data $u_p \in \overline{\shM}_{X,f(p)}^\vee
= \Hom(\overline{\shM}_{X,f(p)}, \NN)$ and 
$u_q \in \Hom(\overline{\shM}_{X,f(q)},\ZZ)$ at marked points $p$ and nodes $q$ of $C$.
\end{itemize}

\begin{definition}
Given a combinatorial type $\tau$ of a stable log map
$f \colon C/W \rightarrow X/S$, we define the associated \textbf{basic monoid} $Q$ by first defining its dual
\begin{equation*}
Q_{\tau}^{\vee}=\left\{ \big((V_{\eta})_{\eta}, (e_q)_q\big)
\in \bigoplus_{\eta} \overline{\shM}_{X,f(\eta)}^{\vee}\oplus \bigoplus_q\NN
\,\bigg|\, \forall q: V_{\eta_2}-V_{\eta_1}=e_qu_q\right\}
\end{equation*}
where the sum is over generic points $\eta$ of $C$ and nodes $q$ of $C$. We then set
\[
Q_{\tau}:=\Hom(Q_{\tau}^{\vee},\NN).
\]
\end{definition}

 Note that, by definition $Q_{\tau}$ indeed only depends on the combinatorial type
of $f \colon C/W \rightarrow X/S$. Geometrically, one can interpret $Q^{\vee}_{\tau,\RR}:=\Hom(Q_\tau,\RR_{\geq 0})$ as the moduli
cone of tropical curves of fixed combinatorial type.

Given a stable log map 
$f \colon C/W \rightarrow X/S$, one can show that there is a canonical map 
$Q \rightarrow \overline{\shM}_W$, where $Q$ is the basic monoid defined by the combinatorial type of $f$.

\begin{definition}
A stable log map 
$f \colon C/W \rightarrow X/S$
is said to be \textbf{basic} if the natural map of monoids 
$Q \rightarrow \overline{M}_W$ is an isomorphism.
\end{definition}

For every $g \in \ZZ_{\geq 0}$, 
$\beta \in H_2(X,\ZZ)$ and $u=(u_1, \dots, u_k)$ with 
$u_i \in |\Sigma(X)|$, we denote by 
$\M_{g,u}(X/S,\beta)$ the moduli space of genus $g$ basic stable log maps to $X/S$ of class $\beta$ and with $k$ marked points of contact data 
\[u=u_1, \dots, u_k.\]

\begin{theorem}[Abramovich--Chen \cite{AbramovichChen}, Gross--Siebert \cite{GSlogGW}, 2011]
 If $X/S$ is projective, then the moduli space 
    $\M_{g,u}(X/S,\beta)$ is a proper Deligne-Mumford stack.
Furthermore, If $X/S$ is log smooth, then the moduli space 
    $\M_{g,u}(X/S,\beta)$ admits a natural virtual fundamental class $[\M_{g,u}(X/S,\beta)]^{virt}$.
\end{theorem}

\subsection{Flow trees to tropical and log curves}
Recall from \S \ref{sec_flow_tree} that quiver DT invariants can be calculated using a universal formula from simpler attractor DT invariants. Further, the coefficients $F_{r,T}^\theta(\gamma_1,\dots,\gamma_r)$ are indexed by attractor tropical trees $T$ in $M_\RR$. In this section we explain how to interpret these coefficients in terms of counts of log curves.

In \cite{NS}, Nishinou-Siebert used toric degenerations and log deformation theory to prove 
a correspondence theorem between counts of genus $0$ complex curves in $n$-dimensional toric varieties and counts of tropical curves in $\RR^n$. This tropical-log correspondence theorem is established in the context when the basic monoid is trivial. In our paper with Bousseau \cite{AB2}, we generalize this tropical-log correspondence to the situation of non-trivial basic monoids. In other words, we work with ``families'' of tropical trees corresponding to log curves over basic monoids of rank equal to the dimension of the family.

Fix a quiver $Q$, a dimension vector $\gamma \in N$ and a decomposition $\gamma=\gamma_1 +\dots+\gamma_r$, such that $\langle \gamma_i,-\rangle \neq 0$ for all $1\leq i \leq r$ and $\langle \gamma,-\rangle \neq 0$.
We denote by $d$ the number of vertices of $Q$, so that $N \simeq \ZZ^d$ and $M_\RR \simeq \RR^d$. Let $X_\Sigma$ be a toric variety as in \S \ref{sec_enum}, with toric divisors $D_i$ containing hypersurfaces $H_i$. Writing $\mathbf{H}=(H_1,\dots,H_r)$, let $\shM_{\mathbf{H}}^{\mathrm{log}}(X_\Sigma)$ be the moduli space of genus $0$ basic stable log maps
 $f: (C, \{ p_1, \ldots , p_{r+1} \} ) \to X_{\Sigma}$ satisfying
\begin{itemize}
    \item[(i)] $f(p_i) \in H_i$ for all $1 \leq i \leq r$ 
  \item[(ii)] The contact order of $f$ with $D_i$ at $p_i$ is  $\langle \gamma_i , - \rangle  \in M$.
  \item[(iii)] The contact order of $f$ with the toric boundary divisor of $X_\Sigma$ at $p_{r+1}$ is $-\langle \gamma, -\rangle$.
  \end{itemize}
  
For general $\mathbf{H}$, one shows in \cite[Theorem 2.8]{AB2} that the moduli space $\shM_{\mathbf{H}}^{\mathrm{log}}(X_\Sigma)$ is $(d-2)$-dimensional and naturally stratified by the combinatorial types of stable log curves. For every combinatorial type $\rho$, there is a stratum $\shM_{\rho,\mathbf{H}}^{\mathrm{log}}(X_\Sigma)$ consisting of stable log maps of type $\rho$, and which is of codimension $\mathrm{rk}\, Q_{\rho}^{\mathrm{gp}}$ in $\shM_{\mathbf{H}}^{\mathrm{log}}(X_\Sigma)$, where  $Q_{\rho}$ is the basic monoid of $\rho$. In particular, types $\rho$ with $\mathrm{rk}\, Q_{\rho}^{\mathrm{gp}}=(d-2)$ correspond to $0$-dimensional strata of $\shM_{\mathbf{H}}^{\mathrm{log}}(X_\Sigma)$, and we denote by $N_{\rho,H}^{\mathrm{toric}}(X_\Sigma)$ the log Gromov--Witten counts of the stable log maps of type $\rho$. The tropicalizations of these stable log maps are $(d-2)$-dimensional families of tropical curves parametrized by the cones $\Hom(Q_\rho,\RR_{\geq 0})$.

\begin{figure}[ht]
\center{\scalebox{.45}{\input{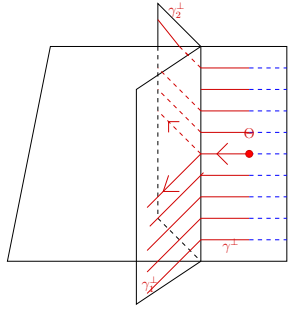}}}
\caption{A family of tropical curves associated to an attractor tree}
\label{families}
\end{figure}

We construct such $(d-2)$-dimensional families $\rho_T$ of tropical curves in $M_\RR$ from an attractor tree $T$, contributing in the formula \eqref{Eq: KS attractor} for quiver DT invariants. To do this, we extend the root of $T$ to infinity to obtain a tropical curve with leaves constrained to lie in the hyperplanes $\gamma_i^{\perp}$. We then consider a deformation of this tropical curve while preserving the combinatorial type and the constraints on the leaves -- see Figure \ref{families}.

We denote by $N_{\rho_T,\mathbf{H}}^{\mathrm{toric}}(X_\Sigma)$ the corresponding log Gromov--Witten invariant counting the finitely many  stable log maps in $\shM_{\mathbf{H}}^{\mathrm{log}}(X_\Sigma)$ of type $\rho_T$-- see \cite[Defn. 1.20]{AB2}. Our main result is \cite[Theorem B]{AB2}:

\begin{theorem} \label{thm}
The coefficients $F_{r,T}^\theta(\gamma_1,\dots,\gamma_r)$ expressing the contribution to $F_r^\theta(\gamma_1,\dots,\gamma_r)$ of an attractor tree $T$ in \ref{Eq: KS attractor} satisfy
\[ F_{r,T}^\theta(\gamma_1,\dots,\gamma_r) = N_{\rho_T,\mathbf{H}}^{\mathrm{toric}}(X_\Sigma)\,.\]
\end{theorem} 

We provide a brief summary of the main steps of the proof here. The first step is to construct a toric degeneration $\mathcal{X} \rightarrow \A^1$ of $X_\Sigma$ and of the constraints $\mathbf{H}$. We then use a degeneration formula expressing the invariants $N_{\rho_T,H}^{\mathrm{toric}}(X_\Sigma)$ of the general fibers $X_\Sigma$ as a sum of invariants $N_{\rho_S}^{\mathrm{toric}}(\mathcal{X}_0)$ of the special fiber $\mathcal{X}_0$, where $S$ are binary trees in $M_\RR$ deforming $T$. We then prove an explicit formula
     \[ N_{\rho_S}^{\mathrm{toric}}(\mathcal{X}_0)=\prod_{v} |\langle \gamma_{v'},
     \gamma_{v''}\rangle | \]
computing $N_{\rho_S}^{\mathrm{toric}}(\mathcal{X}_0)$ as a product over local contributions of the trivalent vertices $v$ of $S$. This formula is proved using the theory of punctured log maps due to Abramovich--Chen--Gross--Siebert \cite{ACGS} to count log maps by gluing punctured log maps. Finally, by the flow tree formula proved in \cite{AB} and briefly reviewed at the end of \S\ref{sec_flow_tree}, one also calculates, 
     \[ F_{r,T}^\theta(\gamma_1,\dots,\gamma_r) = \sum_S \prod_{v} |\langle \gamma_{v'},
     \gamma_{v''}\rangle|\,.\]
We then conclude that $F_{r,T}^\theta(\gamma_1,\dots,\gamma_r) = N_{\rho_T,H}^{\mathrm{toric}}(X_\Sigma)$.

Another way to construct a geometrical object out of a quiver is based on the theory of cluster algebras
\cite{FZ} and cluster varieties \cite{FG}.
Starting with a quiver $Q$ with $d$ vertices, one can construct a $d$-dimensional cluster variety using cluster transformations in a way prescribed by the combinatorics of $Q$. 

In a paper in preparation with Bousseau, we provide a correspondence between quiver DT invariants, in the situation when attractor DT invariants are trivial, and counts of log curves in cluster varieties.
The proof will use a comparison between the
stability scattering diagram of Bridgeland \cite{MR3710055} for quiver DT invariants and the canonical scattering diagram of Gross-Siebert \cite{gross2021canonical} for log Gromov--Witten invariants of log Calabi-Yau varieties, using as intermediate steps the cluster scattering diagram of Gross--Hacking--Keel--Kontsevich \cite{GHKK} and the comparison results obtained in  \cite{HDTV, ABclustermirror}.
This will generalize the Kronecker-Gromov--Witten correspondence previously shown in the $d=2$ case \cite{bousseau2020quantum, GP, GPS,MR3033514, MR3004575}. For variants of the Kronecker-Gromov--Witten correspondence in dimension two see also \cite{bousseau2018example, bousseau2023all, reineke2021moduli}.
Such a correspondence between quiver DT invariants and log curve counts in cluster varieties is compatible with Theorem \ref{thm} using that cluster varieties are obtained by non-toric blow-ups of toric varieties \cite{GHKbirational}, and so admit  degenerations into a toric variety as in \cite{HDTV,GPS}. In particular, log curve counts in cluster varieties are naturally related to log curve counts in toric varieties appearing in Theorem \ref{thm}. 

\bibliographystyle{plain}  
\bibliography{biblio.bib}

\end{document}